\documentclass{article}%[eqsection]

\newenvironment{pf}{\noindent{\sc Proof}.\enspace}{\rule{2mm}{2mm}\medskip}

\newtheorem{theorem}{Theorem}[section]

\newtheorem{lemma}{Lemma}[section]

\newtheorem{remark}{Remark}[section]
\newtheorem{remarks}{Remark}[section]
\newtheorem{definition}{Definition}[section]
\newcommand{\be}{\begin{equation}}
\newcommand{\ee}{\end{equation}}

\newcommand{\N}{{\bf N}}

\renewcommand{\l }{\lambda }

\begin{document}

\title{\bf In whose mind is Mathematics \\an ``\emph{a priori} cognition''?}

\author{Massimiliano Berti $^1$, Antoine Suarez $^2$ and Rocco Tarchini $^3$\\{\it\small $^1$ Dipartimento di Matematica e Applicazioni, Universit\`a di Napoli, Italy, {\tt m.berti@unina.it}}\\ {\it\small $^2$ Center for Quantum
Philosophy, Zurich, Switzerland, {\tt suarez@leman.ch}}\\ {\it\small $^3$ Residenza Universitaria Monterone, Napoli, Italy, {\tt roccotarchini@tiscali.it}}}

\date{September 22, 2008}

\maketitle

\noindent
{\bf Abstract:} According to the philosopher Kant, Mathematics is an ``a priori cognition''. Kant's assumption, together with the unsolvability of Hilbert's $X$-th problem, implies an astonishing result.

\section{Introduction}

\large

In a celebrated contribution on ``Mathematical problems'', at the Second International Mathematicians Congress (Paris, August 8, 1900), David Hilbert shares the conviction that there is an answer to every mathematical problem, and that we can find it by a finite number of purely logical steps:

\begin{description}
\item \quad  ``Take any definite unsolved problem such as [...] the existence of an infinite number of prime numbers of the form $ 2^n + 1 $. However unapproachable these problems may seem to us and however helpless we stand before them, we have nevertheless the firm conviction that their solutions must follow by a finite number of purely logical processes. [...] We hear within us the perpetual call: There is the problem. Seek its solution. You can find it by pure reason, for in mathematics there is no {\it ignorabimus}."\cite{H}
\end{description}

At the Congress Hilbert presented 23 problems. The $X$-th in the list was the following:

\begin{description}
\item	
\qquad {\bf ($X$-th Hilbert problem)} ``Given a diophantine equation with any number of unknowns, devise a process which could determine by a finite number of operations whether the equation is solvable in integers."
\end{description}

Assuming there is an answer to every ``definite" mathematical problem, we are led to the question:
will we have in the future a universal method which allows us to solve any
mathematical problem simply by calculation? Will the day come in which the solution to a given problem becomes only a matter of computing time? This question is the so called \emph{Entscheidungsproblem} \cite{HA}.

\smallskip

Work by K. G\"{o}del (1931), A. Turing (1936-7), A. Church (1935-6), E. Post (1936) and others
(all these basic papers are collected in  \cite{D}),
led to the astonishing result that the answer to Hilbert's \emph{Entscheidungsproblem} is ``no".
We will never have a universal computer program capable of answering any well-posed question on numbers.
This result is usually referred to as an \emph{Undecidability} theorem.

\smallskip

Even more astonishing is that (against Hilbert's conviction) the answer to the $X$-th Hilbert's problem is negative: there is no process for determining by a finite number of operations whether any diophantine equation is solvable. This result, proved in Theorem
\ref{inco2}, is a consequence of the abstract
\emph{Undecidability} theorem  \ref{inco} plus the famous DPRM-theorem \ref{matia} in number theory (DPRM is the acronym for the authors: Davis, Putnam, Robinson, Matiyasevich).

\smallskip

By contrast to the general \emph{Undecidability} problem, $X$-th Hilbert's problem does not refer to abstract logical operations, but involves a very intuitive and sharp defined arithmetical context. This makes it quite appropriate for discussing claims about the status of Mathematics.

In particular we discuss in this article a celebrated claim of the philosopher Immanuel Kant. In the preface to the second edition (1787) of his work \emph{Critique of pure reason}, Kant states that Mathematics is an ``a priori cognition'' (``Erkenntnis a priori''). The mathematician, when he proves a theorem regarding some mathematical object (for instance a triangle), he must not attribute to the object any other properties than those which necessarily follow from the properties he himself has placed into the object, in accordance with his conception.

\smallskip

The \emph{Sections \ref{2}-\ref{4}}
present the main steps that led to the result of the unsolvability of the $X$-th Hilbert's problem. In \emph{Section \ref{5}} we discuss
the implications of this result for Kant's conception of Mathematics as an  ``a priori cognition''.

The interested reader can find in the Appendix a more detailed
presentation of the \emph{Halting problem}, which is crucial
for proving the unsolvability of $X$-th Hilbert's problem.

\section{Hilbert's $X$-th problem and the notion of algorithm}\label{2}

A diophantine equation is an equation of the type $ p(x_1, \ldots, x_k)  = 0 $
where $ p $ is a polynomial with integer coefficients and the solutions $ (x_1, \ldots , x_k) $
are searched in the set of the natural numbers $ \N $. So, for example
$$
2x_1 - 3 x_2 = 1
$$
with $ x_1$, $ x_2 $ in $ \N $
is a diophantine equation  generated by a polynomial of degree one. It is easy to check that all its solutions
are given by $ x_1 = 2 + 3k $ , $ x_2 = 1+ 2 k  $, for all $k $ in $ \N $.

The diophantine equation
$$
x_1^2 - x_2^2 = 3 \, ,
$$
whose unique solution is $ x_1 = 2 $, $ x_2 = 1 $,
is generated by a polynomial of degree two.

Diophantine equations arise in many mathematical contexts.

\smallskip

Any diophantine equation of degree one can be solved by the so called Euclidean algorithm taught at school.
This is a ``mechanical" procedure.

A creative effort of mathematicians has allowed to discover a general procedure to solve also any diophantine equation of degree two.

What about equations of higher degree? As already said, the answer to this question is arduous and involves many advanced tools of number theory.

The first thing to do in order to answer the $X$-th Hilbert problem is to understand what Hilbert meant by ``purely logical processes".

\smallskip

In the last decades of the nineteenth century,
Frege and Peano had showed that
the ``ordinary reasoning" used for the proofs of mathematical theorems amount to a finite number of
formal manipulations of a finite number of symbols, like for example
\begin{description}
\item  quantifiers: $ \forall $  (for all), $ \exists $  (there exists)
\item symbols: $ \in $  (belongs to),  $ \notin $  (does not belong to), $ \subseteq $ (included), $ \subset $ (strictly included),

\end{description}
acting on a finite alphabet, like
the characters 0-9, and letters $ x, n, \ldots $ denoting variables.

This attitude was the so called formalization of arithmetics.

\smallskip

By ``purely logical processes" Hilbert meant exactly these formal manipulations of symbols.
Such %``formal" or ``mechanical''
operations were also referred to as ``deterministic", ``computable", ``effective". The main feature
is that they are expressed
by a finite set of instructions of finite size, and, as such,
they came to be defined as ``algorithmic procedures" or simply ``algorithms".
Several definitions attempting to grasp the essence of the notion of
``algorithm"
have been proposed by
G\"odel-Herbrand, Kleene, Church, and Turing, see e.g. \cite{R}, \cite{Kl}, \cite{T}.
Happily it has been proved that all such definitions are equivalent.
This is the strongest argument to support the famous ``Church's thesis" according to which
these formal definitions translate the intuitive notion of ``algorithm".
In the following we accept the Church's Thesis as granted.

In the Appendix we shortly (but rigorously)
present the Turing characterization which has the advantage to be the most intuitive,
especially nowadays that we are used to standard computer machines.

\smallskip

For the reader who does not wish to learn the mathematical details, it is sufficient to
consider a Turing Machine (TM) as a standard computer program
acting on an idealized computer with an infinite memory storage, working with a finite number
of digits, and which never breaks down.

Given an input value $ x \in \N $, the computer starts performing its instructions.
It can happen that such computations stop giving an output value $ y \in \N $. However, it
could also happen that the computer program does not stop when initialized with a certain input value $ x$.
For example,
it is familiar that sometimes a computer program  ``enters in loop", repeating periodically the very same
instructions. Thus, any TM defines a function $ y := f(x) $, which may not be defined for all $ x \in \N $.
The functions defined by TMs are called  ``general recursive functions".

\smallskip

It is fundamental for the understanding of the following to realize that the class of all the general recursive functions is a countable infinity, i.e. it can be put in one-to-one
correspondence with the natural numbers. The reason is that all the possible computer programs which can be written with a finite number of instructions of finite size form a countable infinity. Moreover the class of the  ``general recursive functions" can be ``effectively" enumerated  as a list
$$
f_n (\cdot ) \, ,  \qquad n = 1, 2, \ldots \, ,
$$
in fact, there is an ``algorithmic" manner to write down such a list (see lemma \ref{enu} below).

\section{The ``Halting problem"}\label{3}

This problem can be formulated as follows:
\begin{itemize}
\item Is there an ``effective" procedure such that, given $ n $ and $ x $ we can
determine whether $ f_n (x)$ is  defined?  In other words, is there an algorithmic way to decide
if the $n$-th Turing machine,  applied to an input $ x $, stops, giving an output value $ f_n (x) $?
\end{itemize}

We say that the  $n$-th Turing Machine\footnote{We identify
the  $n$-th TM with the function $ f_n $ that it produces.} applied to the input value
$ x \in \N $ is convergent if it stops, giving an output value $ f_n (x) $ (``halting" means
``having an output"). On the contrary we say that TM is divergent.

\begin{theorem}\label{thm:halting}
There is no general recursive function (of $2$ variables) such that
\be\label{gx}
g(n,x) := \cases{1 \qquad  {\rm if } \ f_n ( x ) \ {\rm is \ convergent} \cr
0 \qquad   {\rm if } \ f_n ( x ) \ {\rm is \ divergent\, . }}
\ee
\end{theorem}

\begin{pf}
Suppose by contradiction that such a function $ g $ does exist.
Then we can define the new function
\be\label{psin}
\Psi (n) :=  \cases{1 \qquad \qquad \qquad \  {\rm if } \ g(n,n)=0 \cr
{\rm divergent} \ \ \, \qquad   {\rm if } \ g(n,n)=1 \  . }
\ee
The function $ \Psi $ is a general recursive function (of one variable) because
$ g(n,n) \in \{0,1\} $ is algorithmically defined for every $ n $.
Being the function $ \Psi $ a general recursive function
there exists $ n_0 \in \N $
such that
$$
\Psi ( \cdot ) = f_{n_0}( \cdot ) \, .
$$
Then $ \Psi (n_0) = f_{n_0}(n_0) $ which, by (\ref{psin}), is convergent (with value $1$)
if and only if $ g(n_0,n_0) = 0 $.
On the other hand, by (\ref{gx}), $ g(n_0,n_0) = 0 $ if and only if $ f_{n_0} (n_0 )$ is divergent.
This contradiction  proves the Theorem.
\end{pf}

Let us understand properly the previous result.
Any pair $(n,x)$ belongs either to the list of the pairs such that the $n$-th TM applied to the input $ x $ stops, or to the   list of the pairs $(n,x)$ such that the $n$-th TM applied to the input $ x $ does not stop. The non existence theorem \ref{thm:halting} states  the impossibility of writing down effectively such lists via an algorithmic procedure.
In other words: for any pair $(n,x)$ an omniscient mind (God) surely knows and can decide to which of the two lists the pair belongs; but no machine can decide it.

This means that the answer to the ``Halting problem'' is negative: there is no ``effective" procedure such that, given $ n $ and $ x $ we can determine whether $ f_n (x)$ is  defined.

\begin{definition}
A set $ S \subset \N $ is recursively enumerable if it is the range of
a general recursive function. In other words if  $ S $ is the set of the outputs of a Turing Machine.
\end{definition}

The Halting Theorem implies the following important consequence (proved in section \ref{64} of the Appendix).

\begin{theorem}\label{thm:rec}
There exists a set $ K \subset \N $ which is recursively enumerable but whose complementary set $ K^c $ is not
recursively enumerable.
\end{theorem}

\section{G\"odel incompleteness and $X$-th Hilbert problem}\label{4}

A formal system  ${\cal F }$ can be thought as a finite alphabet of symbols,  a finite set of axioms and of
rules of inference, through which it is possible to deduce new propositions,
see \cite{PE} for a wider discussion. We shall  consider ``consistent" formal systems,
in which it is not possible to deduce simultaneously both a proposition and its negation (it is
not possible to prove contradictions from the axioms). We shall consider systems in which
each formally provable proposition is actually true. We shall say that a formal system is
``complete" if every true proposition of the system is  formally provable.

\smallskip

Computability theory provides a perspective from which it can be seen that incompleteness
is a pervasive fundamental property of formal systems. Results of this type were
originally achieved by G\"odel  \cite{G} in 1930.
For definiteness (with no loss of generality) we now consider the formal system of Peano Arithmetics (PA) (it is nothing but
a formalization of the ordinary arithmetics, see \cite{C})  and we prove
that it is not complete.

\smallskip

Let $ K $ be the  set defined in theorem \ref{thm:rec}.
The statement $ n \notin K $ is a proposition of the (PA)-formal system, say $ P_n $.

\begin{theorem} {\bf (Incompleteness Theorem)} \label{inco}
There is $ n_0 $ such that the statement
$$
n_0 \notin K \quad {\rm is \ true}
$$
but $P_{n_0} $ is  not provable in the arithmetic formal system (PA).
\end{theorem}

\begin{pf}
By contradiction suppose that, for all $  n $,
\be\label{sopra}
n \notin K \quad {\rm is \ true}  \quad {\rm if \ and \ only \ if } \quad P_n {\rm \ is \ provable \  in \ the \ formal \ system}\, .
\ee
We now prove that the complementary set $ K^c $ is recursively enumerable, contradicting
theorem \ref{thm:rec}.
For every $ n \in K^c $, the statement $ n \notin K $ is true, and,  by (\ref{sopra}),
the proposition $ P_n $ is formally provable.
But the set of all the formally provable propositions is recursively enumerable. Indeed
we can algorithmically list all the  formal propositions (for example
using a lexicographical ordering). Next we can algorithmically check whether
each formal proposition obeys  the formal rules of inference in PA. Then we
have a purely mechanical test to know  if each formal proposition is
a theorem.
\end{pf}

We can now consider  the original $ X $-th Hilbert problem. After Turing's achievements about the halting problem,
the complete solution of the unsolvability of the $ X $-th Hilbert problem required  a long time.
The main breakthrough was the following theorem due to joint
work of J. Robinson, M. Davis, H. Putnam, see \cite{DMR}, and Y. Matiyasevich,\cite{YM}.

\begin{theorem}  {\bf (\cite{DMR},\cite{YM})}\label{matia}
Let $ S \subset \N $ be a recursively enumerable set. Then there exists
a polynomial $ p (n, x_1, \ldots , x_k) $ with integer coefficients such that
the equation
\be\label{pne}
p (n, x_1, \ldots , x_k) = 0  \ \  { \rm has \ integer \ solutions} \quad {\rm if \ and \ only \ if }  \quad n \in S \, .
\ee
\end{theorem}

We underline that the construction of the polynomial $ p $ is completely explicit.
Let $ p_0 $ be the polynomial provided by theorem \ref{matia} when $ S = K $ is the set
defined in theorem \ref{thm:rec}. By theorem \ref{inco} we deduce the following G\"odel type
result in arithmetics.

\begin{theorem} {\bf (Diophantine Incompleteness Theorem)} \label{inco2}
There is $ n_0 $
such that the equation
$$
p_0 (n_0,x_1, \ldots, x_k) = 0   \ \  { \rm has \ no \ integer \ solutions}
$$
although such a statement is not formally provable in the arithmetics.
\end{theorem}

Finally, let us conclude proving the negative answer to the $X$-th Hilbert problem.

\begin{theorem}\label{noH}
There exists a class of diophantine equations for which it is not possible to decide algorithmically
whether they admit integer solutions.
\end{theorem}

\begin{pf}
The set $ K^c $ is not recursively enumerable by theorem \ref{thm:rec}.
Then, given $ n_0 \notin K $,
 it is not possible to determine algorithmically that
$ n_0 $ actually  is not in $ K $.
Therefore, by (\ref{pne}), we cannot determine
algorithmically  if the equation $ p_0 (n_0,x_1, \ldots, x_k) = 0 $ has  or  not  integer solutions.
\end{pf}

\section{Philosophical conclusions}\label{5}

To discuss Kant's claim on the status of Mathematics, we point out the following consequence of the previous results: %Theorem \ref{noH}:

\begin{description}
\item \quad \ ``At any time $ T $ there is a diophantine equation $ E $ such that no human being at time $ T $
can be sure to be able to answer the question whether $ E $ is solvable or not".
\end{description}

Indeed, at any time $ T $, there exists no algorithm to generate $ K^c $ (Theorem \ref{thm:rec}).
Therefore,  for any given  $ n_0 \in K^c $,
we shall never possess an algorithm which is able to give us the answer: $``n_0 \in K^c"$.

After some mental work, at some time $ T' $ after  $  T $,   we  might perhaps get the ``insight"  for developing a more powerful procedure allowing us to state that $``n_0 \in K^c"$. But at time $ T $ we cannot be sure to receive such an insight! Therefore, at any time $ T $ there exists $ n_0 \in K^c $, such that I'm not sure to be able to answer: ``the diophantine equation (E)
$ p_0(n_0, x_1, \ldots, x_k ) = 0 $ has no solutions".

\smallskip

In more simple words: at any time $ T $ the question whether any diophantine equation  is solvable or not, is in principle answerable and the answer exists: it is either YES or NO. But, for the previous diophantine equation $ E $, the answer
doesn't exist in any human mind, not even in the implicit form of a result like:
\emph{``I can be sure} to obtain it after a certain finite number of operations of finite size".

\smallskip

What are the consequences of this conclusion with relation to Kant's conception of Mathematics as an ``a priori'' cognition?
Kant presents his conception mainly in the two following quotations:

\begin{description}
\item \quad  ``A new light must have flashed on the mind of the first man (Thales, or whatever may have been his name) who demonstrated the properties of the isosceles triangle. For he found that it was not sufficient to meditate on the figure, as it lay before his eyes, or the conception of it, as it existed in his mind [...] but that it was necessary to produce these properties, as it were, by a positive a priori construction; and that, in order to arrive with certainty at a priori cognition, he must not attribute to the object any other properties than those which necessarily followed from that which he had himself, in accordance with his conception, placed in the object.  [...]

    A much longer period elapsed before physics entered on the highway of science [...] When Galilei experimented with balls of a definite weight on the inclined plane, when Torricelli caused the air to sustain a weight which he had calculated beforehand to be equal to that of a definite column of water [...] a light broke upon all natural philosophers. They learned that reason only perceives that which it produces after its own design;"    \cite{K}
\end{description}

Thus, according to Kant Mathematics and the mathematical description of Nature is something existing ``a priori'' in some mind independently of experience. Additionally, he seems to suggest that this mind is a human one (the mind of Thales, Galilei, Torricelli, etc.).

However, the fact that at any time T there is a diophantine equation E such that no human can say whether it is solvable or not, even if the answer to this question exists, proves that mathematical truth will never be completely contained in any human mind. Therefore Mathematics does not exist  ``a priori'' in a human mind.

\smallskip

Then, if we maintain Kant's conception that Mathematics is an ``a priori'' cognition, we are led to conclude that it is ``a priori'' in some other mind, who is mightier than the human one, a mind who actually contains the whole of Mathematics at once.

In this sense one can say that the existence of an omniscient mind (God) follows from Kant's assumption on Mathematics in his \emph{Critique of pure reason} and the recent theorems on the foundations of Mathematics.

\smallskip

\newpage

\section{Appendix}

For the reader convenience we report the basic ideas on ``general recursive" functions.

\subsection{The primitive recursive functions}

The following class of functions, that we shall  call the
``primitive recursive functions",  must be considered  ``algorithmic", ``effectively computable".

\begin{itemize}
\item ($i$) The constant functions $f(x_1, \ldots , x_k) := c $.
\item ($ii$) The successor function $ f(x) := x +1 $.
\item ($iii$) The identity functions $ f(x_1, \ldots , x_k) := x_i $,  $ i = 1, \ldots, k $.
\end{itemize}

Then we can produce new recursive functions via ``composition" (composition of primitive recursive functions
yields a new one)
and ``induction".

Deterministic operations like ``sum" and ``multiplication" are primitive recursive functions
as well as every polynomial. However, it is possible to construct functions,
intuitively computable by an algorithm,
which are not primitive recursive, see e.g. \cite{R}.
The class of primitive recursive functions
does not exhaust the concept of  ``effectively computable" functions.

\subsection{The general recursive functions}

At first glance an insuperable obstacle to give
a definition of  ``effectively computable" functions (that we shall call
 ``general recursive" functions) appears.
Indeed, suppose one had an  ``effective" definition of general recursive functions.  Since there are only
countably many instructions one can give for performing a calculation, one could enumerate
these functions, say $ f_n (x) $, $ n =1, 2, \ldots $. Then we can define the new function
\be\label{gn}
g(n) := f_n(n) +1 \, .
\ee
Such a function $ g $ is clearly produced by an algorithmic procedure. Hence $g$ belongs to the previous list
of general recursive functions, i.e.
\be\label{gf}
\exists \, n_0 \ \  {\rm  such \ that} \quad
g (\cdot ) = f_{n_0} (\cdot ) \, .
\ee
But, let us consider the value $ g( n_0 ) $. On one side, by (\ref{gf}),
$ g ( n_0 ) = f_{n_0} ( n_0 ) $. On the other side, by (\ref{gn}),
$ g ( n_0 ) = f_{n_0} ( n_0 ) + 1 $. Such a contradiction proves that an ``effective" definiton
of general recursive functions does not exist.

%\begin{remark}
%The previous procedure is an example of the ``Cantor diagonal argument".
%\end{remark}
\smallskip

The way out of this dilemma is to give a definition of ``general recursive" functions which
are not defined for all the values $ x \in \N $.  We shall call such functions ``partial  recursive" functions.
\\[1mm]
{\bf Example}. Consider the function $ f(x)  $  defined as follows: carry out the decimal expansion of $ \pi $ until
a run of at least $ x $ consecutive $5$ appears; if and when this occurs
define $ f(x)$ as the position of the first digit of this run.
We do not know, at the moment,
whether this algorithmic procedure produces a function defined for all $ x \in \N $.
Such $ f $ should be considered a partial recursive function.

\smallskip

For partial general recursive functions,
the previous argument does not lead to a contradiction because
$ g( n_0 )$ does not need to have a value, namely $ n_0 $ could not belong to the
 domain of definition of the function $ g $.

 \smallskip

We present below the Turing characterization of general recursive functions.

\subsection{The Turing machines}

A Turing machine is an idealized universal computer device
which can operate without restrictions of memory storage and space.
It can be thought as a box, containing a finite set of internal states
$ q_0,  q_1,  \ldots q_k $,
together  with an infinite tape.  The Turing machine (TM) can perform the following basic operations
which depend on the internal states and on the input value read in the cell tape being examined:
\begin{description}
\item (i) A TM can write $ ``1" $ on the cell of the tape it is examining  if there is no yet $ ``1" $ printed.
\item (ii) A TM can erase $ ``1" $ on the cell of the tape it is examining  if there is $ ``1" $ printed.
\item (iii) A TM can shift its attention one cell to the right, or to the left,  along the tape.
\end{description}

More precisely  a TM can be defined as a finite set of
ordered quadruples consisting of symbols for:
\begin{description}
\item (i) an internal state
\item (ii) a possible tape cell input value
\item (iii) an operation
\item (iv) an internal state.
\end{description}

Such a  quadruple  expresses an instruction for a TM:
given the internal state (i) and
the tape cell input value (ii), the machine performs the operation (iii)
and takes the internal state (iv) (this set of instructions  can be thought as a computer program).

A TM is a deterministic device, that is, given an internal state (i) and a tape cell input value (ii), the
operation performed in (iii) is uniquely determined (there are not $2$ different quadruples
having the same two symbols (i) and (ii)).

\smallskip

Now we can define the following class of functions $ y = f(x) $ of one integer variable.
Given $ x \in \N $, build a tape with $ x $ consecutive $1$ digits (and $ 0 $ outside).
This tape is the input for a specific TM which has the initial internal state $ q_0 $.
The TM  examines first the cell containing the first $ 1$ digit. Then the TM performs its calculus and if and when
it stops take as $ y := f(x) $ the number of digits $1$ occurring in the final tape.
In this case any TM defines uniquely a function
$$
f: {\cal D}(f) \subseteq \N \to \N
$$
which is properly defined only on a  subset $  {\cal D}(f) \subseteq \N $. If, given $ x \in \N $,
the calculus of the TM never stops, we cannot define the value $ f(x)$.
The class of all the functions produced by a TM is called the class
of the ``partial recursive functions".

It is clear that the previous definition of  ``partial recursive functions" can be extended
to functions with more than one integer variable.

\begin{remark}
Recalling the previous example, it is clear that it is possible to write a TM (computer program) to
check if, given any $ x $, there exists a run of $ x $-consecutive $5$ in the decimal expansion
of $ \pi $. We do not know if, given $ x $, such a program stops, providing an answer.
\end{remark}

\begin{remark}\label{rem1}
There are TMs which, for any $ x \in \N $, never stop (in such a case they do not define any function).
An example of such ``inefficient" Turing Machine is provided by $ q_011q_0 $. Indeed the previous set of instructions define
a ``loop": such a machine keeps on reading and printing
indefinitely  the value $ 1 $ on the same cell tape.
\end{remark}

It can be easily recognized that one can define Turing Machines to perform
primitive operations like sum, multiplication, composition, etc... (that is all the primitive recursive functions).
We refer to \cite{R} for an explicit example of a TM machine programmed to calculate
the function $ x \mapsto 2x $. With a bit
of expertise one is led to define a function to be % ``effective",
%``computable",
%``effectively computable", ``recursively computable",
``effectively computable",
%``produced by an algorithm",
if it is defined by a Turing Machine.

\smallskip

\begin{lemma}\label{enu}
The class of  the partial recursive functions is countable and can be
``effectively" enumerated  as a list $ f_n $, $ n = 1, 2, \ldots $.
\end{lemma}

\begin{pf}
The set of Turing Machines
is countable, because it is produced by a finite number of instructions of finite size.
We can arrange all the Turing Machines through
a lexicographical ordering, listing first the TM's defined by $ 1 $ quadruple, next
the TM's defined by $ 2 $ quadruples, by $ 3 $ quadruples, and so on.
Such a procedure is purely algorithmic.
\end{pf}

\subsection{Proof of Theorem \ref{thm:rec}}\label{64}

Define the function
$$
\Phi (n,x) := \cases{(n,x) \qquad \qquad  {\rm if} \ f_n(x) \ {\rm is \ convergent}  \cr
{\rm divergent} \qquad  \  {\rm if} \ f_n(x) \ {\rm \ is \ divergent} \, . }
$$
This function $ \Phi $ is partial recursive.  Indeed, given $(n,x)$ we can find
algorithmically the function $ f_n $ (see lemma \ref{enu})
and compute $ f_n (x)$. If  $ f_n (x)$ is convergent
then we print $ (n,x) $ which is the output of the algorithm $ \Phi $. Otherwise, if $ f_n (x)$ is divergent then also
$ \Phi (n,x) $ is divergent.

Then $ K = {\rm range}(\Phi) $ is recursively enumerable (by definition). If also $ K^c $ were
recursively enumerable then it would be possible to decide algorithmically given any $(n,x)$
if $ f_n (x) $ is convergent or divergent, contradicting
the Halting Theorem \ref{thm:halting}.

\begin{remark}\label{3.1new} By Theorem \ref{thm:rec}
we immediately get  the existence of  a subset $ K \subset \N $ which is recursively enumerable yet not
recursive, by a computable bijection $ \N^2 \leftrightarrow \N $.
\end{remark}

\end{document}